\documentclass{ifacconf}

\usepackage{graphicx}      % include this line if your document contains figures
\usepackage{natbib}
\usepackage{amsmath,amssymb,mathrsfs}
\usepackage{enumerate}
\usepackage{setspace}
\usepackage{bm}
\newtheorem{Theorem}{Theorem}

\begin{document}
\begin{frontmatter}
\title{On the Capability of PID Control for Nonlinear Uncertain Systems}
\author[First]{Cheng ZHAO and Lei GUO}
\address[First]{Institute of Systems Science, AMSS, Chinese Academy of Sciences,
Beijing, 100190 China, (e-mail: 1154279089@qq.com, lguo@amss.ac.cn).}

\begin{abstract}                % Abstract of not more than 250 words.
It is well-known that the classical PID controller is by far the most widely used ones in industrial processes, despite of the remarkable progresses of the modern control theory over the past half a century. It is also true that the existing theoretical studies on PID control mainly focus on linear systems, although most of the practical control systems are inherently nonlinear with uncertainties. Thus, a natural question is: can we establish a theory on PID controller for nonlinear uncertain dynamical systems? This paper will initiate an investigation on this fundamental problem, showing that any second order uncertain nonlinear dynamical systems can be stabilized globally by the PID controller as long as the nonlinearity satisfies a Lipschitz condition. We will also demonstrate that this result can be generalized neither to systems with order higher than $2$, and nor to systems with nonlinear growth rate faster than linear in general.
\end{abstract}


\end{frontmatter}
%===============================================================================

\section{Introduction}
Over the past half a century, remarkable progress has been made in modern control theory and its applications. Despite of this, the classical PID (proportional-integral- derivative) controller (or its minor variations) is recognized to be the most widely used controller in engineering systems by far. For example, in process control, more than 95\% of the control loops are of PID type, and most loops are actually PI control, see \cite{Astrom1995,Astrom2006}.

There are several reasons for the effectiveness of the PID controller: the implementation of the PID controller does not need precise
mathematical models; it can reduce the influence of the  system uncertainties by feedback signals including the proportional action;
it has the ability to eliminate steady state offsets via the integral  action; and it can also  anticipate the future tendency through the derivative action.
Also,  the celebrated Newton's second law in mechanics still plays a fundamental role in modelling dynamical systems of the physical world, which is
actually a second order differential equation of the position of a moving body, and  can be well  regulated by the PID controller, as is well known
in practice,  and as will be justified rigorously in the current  paper.  Of course,  one of the most challenging tasks for the implementation of the PID
control is how to design the three parameters of the controller, which has also been investigated extensively in the literature but most for linear systems, see \cite{Astrom1995, Astrom2006, Blanchini2004, Hara2006, Ho1997, Ho2003, Keel2008, Killingsworth2006, silva2007pid, Soylemez2003}, among others.

One of the most eminent methods for designing the PID parameters is the Ziegler-Nichols rules (see \cite{ziegler1942}), which  is essentially an experimental method. Many other
methods including tuning and adaptation for the design of the PID parameters have also been proposed and investigated, see \cite{Astrom1995,Astrom2006}.
To the best of our knowledge, there is few theoretical results on PID control for nonlinear uncertain dynamical systems concerning global stabilization and
control performance in the literature. However, in order to justify the remarkably practically effectiveness of the PID controller, we have to face with such complex systems, since nonlinearity and uncertainty are ubiquitous for real world systems. This is a longstanding problem in control theory (see, e.g. \cite{Guo2011}) and is the prime motivation of the current paper.

Following a similar theoretical framework  as    the investigation of the maximum capability of the feedback mechanism in \cite{Xie2000,Guo2014}, we will in this paper investigate the capability and limitations of the PID controller  in dealing with  nonlinear uncertain dynamical systems.  We will show that second order nonlinear uncertain dynamical systems will be stabilized globally by PID controller with fixed control parameters, as long as the related nonlinear  uncertain function  satisfies a Lipschitz condition with arbitrary Lipschitz constant.
By doing so, we will also be able to determine an open unbounded set in $\mathbb{R}^3$, from which the stabilizing PID controller parameters can be chosen conveniently.
Moreover, we will further show that such a nice theoretical result cannot be extended to uncertain dynamical systems with  growth rate of the nonlinearity faster than linear,
nor to dynamical systems with order higher than $2$. These somewhat natural results clearly demonstrate the limitations of the classical PID control.

      The rest of the paper is organized as follows. The problem formulation will be described  in the next section. Section $3$ will presents our  main results, with their
proofs put in Section $4$. Section $5$ will concludes the paper with some remarks.

\section{Problem Formulation}
Let us consider a moving body in  $\mathbb{R}^n$ which is regarded as a controlled system. Denote  $x(t)$, $v(t)$, $a(t)$ be its position, velocity, acceleration at the time instant $t$,  respectively. Assume that the external forces acting on the body consist of $f$ and $u$, where $f=f(x,v)$ is a nonlinear function of the position and velocity and $u$ is the control force. There are many examples which satisfy these assumptions. Classical examples contain spring oscillator, pendulum, damped vibration, etc.

By Newton's second law, we have the equation
\begin{equation}
 ma=f(x,v)+u
 \end{equation}
where $u$ is the control input and $m$ is the mass.
Our control objective is to design an output feedback controller to guarantee that for any initial position and initial velocity, the position trajectory tracks a given reference value  $y^{*}\in\mathbb{R}^n$  and at the same time the velocity of the body tends to $0$.

In this paper, our control force is described by the classical PID controller
\begin{equation}
u(t)=k_{p} e(t)+k_{i} \int_{0}^{t} e(s)ds+k_{d} \overset{.}{e}(t),
\end{equation}
where $e$ is the control error $$e(t)=x(t)-y^{*}.$$ The control variable is thus a sum of three terms: the P-term(which is proportional to the error), the I-term(which is proportional to the integral of the error) and the D-term(which is proportional to the derivative of the error). Without loss of generality, we assume that the body has the unit mass $m=1$. Notice that $v=\overset{.}{x}$,  $a=\overset{..}{x}$, then $(1)$ can be rewritten as
$$\overset{..}{x}=f(x,\overset{.}{x})+u.$$
Denote $x_{1}=x$ and $x_{2}=\overset{.}{x}$, the state space equation of this basic mechanic system is
\begin{equation}
\begin{cases}
\overset{.}{x}_{1}=x_{2}\\
\overset{.}{x}_{2}=f(x_{1},x_{2})+u\\
u=k_{p} e(t)+k_{i} \int_{0}^{t} e(s)ds+k_{d} \overset{.}{e}(t)
\end{cases}
\end{equation}
where $x_{1}(0),x_{2}(0)\in \mathbb{R}^{n}$.

 In this paper, we will show that the three controller parameters $k_{p},k_{i},k_{d}$ can be designed such that the position of the body tracks a given $y^{*}$ under the control law $(2)$ for any initial position and velocity, as long as $f=f(x,v)$ is a Lipschitz continuous function with a known Lipschitz constant $L$.

\section{The Main Results}

The performance of the closed-loop system under the  PID controller can be described by the following theorem, established in $n$-dimensional space with $n\ge1$.

Firstly, we define a functional space:
 \begin{small}\begin{align*}\mathscr{F}_{L}=\begin{small}\{f:\mathbb{R}^{2n} \rightarrow \mathbb{R}^n \big| \|f(x)-f(y)\| \leq L\|x-y\|,\forall x,y\in \mathbb{R}^{2n}\}\end{small}\end{align*}\end{small} where $\|.\|$ is the standard Euclidean norm.
\vskip 0.5cm
\begin{Theorem}
Consider the PID controlled system $(3)$ with any unknown $f \in \mathscr{F}_{L}$. Then for any $L>0$, there exists an unbounded open set $\Omega_{K}$ in $\mathbb{R}^{3}$, such that whenever the  controller parameters $(k_{p},k_{i},k_{d})$ are taken from $\Omega_{K}$,  the closed-loop system $(3)$ will be globally stable and satisfies $\lim_{t \to \infty} x_{1}(t)=y^{*}$ and $\lim_{t \to \infty} x_{2}(t)=0$ with an exponential  rate of convergence for any initial value $(x_{1}(0),x_{2}(0)) \in \mathbb{R}^{2n}$, where $y^{*} \in \mathbb{R}^{n}$ is any given setpoint in $\mathbb{R}^{n}$.
\end{Theorem}
\vskip 0.5cm
\begin{rem}
  In fact, the selection of the three controller parameters has wide flexibility and is robust to some extent, due to the open property of the parameter domain $\Omega_{K}$. A small perturbations of these parameters do not change the qualitative performance of the system. From the proof of Theorem $1$, the concrete definition of $\Omega_{K}$    in $\mathbb{R}^{3}$ can be taken as,
\begin{align*}
\Omega_{K}=\bigg\{\begin{bmatrix}k_{p}\\k_{i}\\k_{d}\end{bmatrix}=\begin{bmatrix}-(\lambda_{1}\lambda_{2}+\lambda_{1}\lambda_{3}+\lambda_{2}\lambda_{3})\\\lambda_{1}\lambda_{2}\lambda_{3}\\
\lambda_{1}+\lambda_{2}+\lambda_{3}\end{bmatrix} \bigg| \begin{bmatrix}\lambda_{1}\\
\lambda_{2}\\\lambda_{3}\end{bmatrix}\in
\Omega_{\Lambda}\bigg\}
\end{align*}
with
  $$\Omega_{\Lambda}=\{\Lambda \big | L\phi(\Lambda)h(\Lambda)<1; \lambda_{i}<0,i=1,2,3; \lambda_{i}\neq\lambda_{j},i\neq j\},$$
where $\Lambda$ is defined as $\Lambda=(\lambda_{1},\lambda_{2},\lambda_{3})$ and
 $$h(\Lambda)=\sqrt{3+\lambda_{1}^{2}+\lambda_{2}^{2}+\frac{1}{\lambda_{3}^{2}}},$$
    $$\phi(\Lambda)=\sqrt{ \frac{(\lambda_{3}-\lambda_{2})^{2}+(\lambda_{3}-\lambda_{1})^{2}+\lambda_{3}^{2}(\lambda_{2}-\lambda_{1})^{2} }{(\lambda_{3}-\lambda_{1})^{2}(\lambda_{2}-\lambda_{1})^{2}(\lambda_{3}-\lambda_{2})^{2} } }.$$

  To simplify this parameter set, we next give a corollary which provides a direct formula for determining $k_{p},k_{i}$ and $k_{d}$.
\end{rem}

\begin{cor}
The following 2-dimensional manifold $\Omega_{K}^{'}$ ,
$$\Omega_{K}^{'}=\biggl\{\begin{bmatrix}k_{p}\\k_{i}\\k_{d}\end{bmatrix} \bigg | \begin{cases} k_{p}=-(\epsilon(1+\epsilon)+(1+2\epsilon) a)\\k_{i}=-\epsilon(1+\epsilon)a\quad\quad0<\epsilon<\frac{1}{4} \\k_{d}=-(a+1+2\epsilon)\quad a>\max\{5L,5\}\end{cases} \biggl\}$$
is a subset of $\Omega_{K}$. Hence if we take $(k_{p},k_{i},k_{d})\in \Omega_{K}^{'}$, then under the PID controller $(2)$, we have for any $ f \in \mathscr{F}_L$ and for any $ y^{*} \in \mathbb{R}^n$, the closed-loop system $(3)$ satisfies $\lim_{t \to \infty} x_{1}(t)=y^{*}$ and $\lim_{t \to \infty} x_{2}(t)=0$ with an exponential rate of convergence, for all initial conditions $(x_{1}(0),x_{2}(0))\in \mathbb{R}^{2n}$.
\end{cor}

\begin{rem}
From the above corollary, we find that the integral parameter $k_{i}$ of the PID controller can be taken arbitrarily small, but cannot be zero, since otherwise there will be  no integral action.  At the same time, we notice that one can choose $k_{p},k_{d}$, which is $O(L)$ for large $L$. Of course, the formula given in  Corollary $1$ is not the unique one. In fact, we  have many choices based on Theorem $1$ according to different requirements.
\end{rem}

\begin{rem} It is worth noting that Theorem 1 gives a global convergence result. At the same time, the selection of the three parameters does not rely on the initial values (position and velocity) and the setpoint $y^{*}$. A natural question is:  whether or not the global Lipschitz condition is necessary for global stabilization?  The answer is "yes" in general. In fact, we can show that if $f$ merely satisfies a local Lipschitz condition, for example $f(x)=\|x\|^{1+\epsilon}$, we cannot expect Theorem $1$ holds globally for any $\epsilon>0$. The following proposition rigorously proves this fact.
\end{rem}
\begin{prop}
Consider the  following nonlinear system,
\begin{equation}
\begin{cases}
\overset{.}{x}_{1}={x}_{2}\\
\overset{.}{x}_{2}=\|x\|^{1+\epsilon}+u\\
u={k}_{p} e(t)+k_i \int_{0}^{t} e(s)ds +{k}_{d} \overset{.}{e}(t)
\end{cases}
\end{equation}

where $x=(x_{1},x_{2})\in \mathbb{R}^{2}$, $e=x_{1}-y^{*}$ and $y^{*}$ is a constant. Then for any $ \epsilon > 0$, any $y^{*}\in \mathbb{R}$ and any $ k_{p},k_{i},k_{d}$, there exists $x(0)=(x_{1}(0),x_{2}(0))\in \mathbb{R}^{2}$, such that the maximal existence interval $[0,a)$ of the closed-loop equation $(4)$ is finite  and  satisfies $$e(t)\ge |e(0)|>0$$ for $t\in [0,a)$.
If in addition $\epsilon \le 1$, we have: $$\lim_{t\to a}e(t)=\infty.$$
\end{prop}

The above proposition may not be surprising because we cannot expect to use PID controllers (which is a linear feedback) to solve the global tracking problem of essentially nonlinear systems, even when we know exactly the dynamics of the system.

 In the final part of this section, we show that it is not possible in general for the PID controllers to globally stabilize nonlinear uncertain systems described by differential equations of order $\ge 3$. For simplifying our proof, we only consider the case $n=3$.

\begin{prop}
Define $$\mathscr{G}_{L}=\{f:\mathbb{R}^{3}\rightarrow \mathbb{R} \big | |f(x)-f(y)|\le L\|x-y\|,\forall x,y \in \mathbb{R}^{3}\}.$$
Consider a PID controlled system as follows:
\begin{equation}
\begin{cases}
\overset{.}{x}_{1}&=x_2\\
\overset{.}{x}_2&=x_3\\
\overset{.}{x}_3&=f(x_1,x_2,x_3)+k_p e(t)+k_i\int_{0}^{t} e(s)ds +k_d \overset{.}{e}(t)
\end{cases}
\end{equation}
where $e=x_{1}-y^{*}$ and $y^{*}$ is a constant.
Then, for any $L>0$, any $ k_{p},k_{i},k_{d} \in \mathbb{R}$ and any $ y^{*} \in \mathbb{R}$,
there exists a function $f \in \mathscr{G}_L$ and an initial value $x(0)=(x_{1}(0),x_{2}(0),x_{3}(0))\in \mathbb{R}^{3}$,
such that the closed-loop system  $(5)$ satisfies
$$\sup_{t \ge 0} |e(t)|=\infty.$$
\end{prop}

\section{Proofs of the main results}

\textbf{Proof of Theorem 1:}
First, we introduce some notations. Denote $y_{00}=\int_{0}^{t} e(s)ds$,  $y_{1}(t)=e(t)$, $y_{2}(t) =\overset{.}e (t)$, then $(3)$ is equivalent to
\begin{equation}
\begin{cases}
\overset{.}{y}_{00}&=y_1\\
\overset{.}{y}_{1}&=y_2\\
\overset{.}{y}_{2}&=f(y_1+y^{*},y_2)-f(y^{*},0)+k_{i}(y_{00}+\dfrac{f(y^*,0)}{k_i})\\&+k_{p}y_1+k_{d}y_{2}
\end{cases}
\end{equation}
Denote $g(y_1,y_2)=f(y_{1}+y^{*},y_2)-f(y^{*},0)$, $y_0=y_{00}+\dfrac{f(y^{*},0)}{k_i}$, we can get a more compact equation,
\begin{equation}
\begin{cases}
\overset{.}{y}_{0}&=y_1\\
\overset{.}{y}_1&=y_2\\
\overset{.}{y}_2&=g(y_1,y_2)+k_{i}y_0+k_{p}y_1+k_{d}y_2.
\end{cases}
\end{equation}
Now, by $f\in \mathscr{F}_L $, it is easy to see that $g\in \mathscr{F}_L$ and $g(0)=0$. Hence $0$ is an equilibrium of $(7)$. Denote

$Y=(y_{0}^{T},y_{1}^{T},y_{2}^{T})^{T}$, $Y^{'}=(y_{1}^{T},y_{2}^{T})^{T}$,
$A=\begin{bmatrix}
0&I&0\\
0&0&I\\
k_{i}I&k_{p}I&k_{d}I
\end{bmatrix}$.

Here $A$ is a $3n\times3n$ matrix and $I$ is an $n\times n$ unit matrix. Then $(7)$ can be rewritten as
\begin{equation}
\overset{.}{Y}=AY+\begin{bmatrix}
0\\0\\g(y_1,y_2)
\end{bmatrix}.
\end{equation}
By simple calculations, using the properties of determinants, we have
\begin{align*}
det(\lambda I_{3n*3n}-A)
&=\begin{vmatrix}
\lambda I&-I&0\\
0&\lambda I&-I\\
-k_{i}I&-k_{p}I&(\lambda-k_{d})I
\end{vmatrix}\\
&=\begin{vmatrix}
\lambda I&0&0\\
0&\lambda I&-I\\
-k_{i}I&-(\dfrac{k_{i}}{\lambda}+k_{p})I&(\lambda-k_{d})I
\end{vmatrix}\\
&=\begin{vmatrix}
\lambda I&0&0\\
0&\lambda I&0\\
*&*&(\lambda-k_{d}-\dfrac{k_{p}}{\lambda}-\dfrac{k_{i}}{\lambda^{2}})I
\end{vmatrix}\\
&=(\lambda^{3} -k_{d}\lambda^{2} -k_{p}\lambda -k_{i})^{n}.
\end{align*}
Take $k_{p},k_i,k_d$ such that $\lambda^{3} -k_{d}\lambda^{2} -k_{p}\lambda -k_{i}=0$ has three distinct negative real roots $\lambda_{1}, \lambda_{2}, \lambda_{3}$. This is feasible because we can adjust all the coefficients of the cubic equation.
Define three matrices

$$P
=\begin{bmatrix}
\frac{I}{\lambda_{1}}&\frac{I}{\lambda_{2}}&\frac{1}{\lambda_{3}^2}I\\
I&I&\frac{1}{\lambda_{3}}I\\
\lambda_{1}I&\lambda_{2}I& I
\end{bmatrix},$$

$$P^{'}=
\begin{bmatrix}
I&I&\frac{1}{\lambda_{3}}I\\
\lambda_{1}I&\lambda_{2}I&I
\end{bmatrix},J
=\begin{bmatrix}
\lambda_{1}I&0&0\\
0&\lambda_{2}I&0\\
0&0&\lambda_{3}I
\end{bmatrix},$$

 then it is not difficult to see that

$$P^{-1}
=\begin{bmatrix}
*&*&\dfrac{\lambda_{1}I}{(\lambda_{3}-\lambda_{1})(\lambda_{2}-\lambda_{1})}\\
*&*&\dfrac{\lambda_{2}I}{(\lambda_{3}-\lambda_{2})(\lambda_{1}-\lambda_{2})}\\
*&*&\dfrac{\lambda_{3}^{2}I}{(\lambda_{3}-\lambda_{1})(\lambda_{3}-\lambda_{2})}
\end{bmatrix}$$
where those $*$ in the elements of $P^{-1}$ means that we don't care about what they are in our proof of the theorem.

Define an invertible linear transformation $Y=PZ$, and denote $Z=(z_{0}^{T},z_{1}^{T},z_{2}^{T})^{T}$, where $z_{0},z_{1},z_{2}$ are  $n$-dimensional column vectors. By the relationship $A=PJP^{-1}$, we can write the equation $(8)$ in a diagonal form,
\begin{align*}
\overset{.}{Z}=JZ+P^{-1}\begin{bmatrix}
0\\0\\g(y_1,y_2)
\end{bmatrix}=JZ+P^{-1}\begin{bmatrix}
0\\0\\g(P^{'}Z)
\end{bmatrix}.
\end{align*}

Consequently, we have
\begin{equation}
\begin{cases}
\overset{.}{z}_{0}&=\lambda_{1}z_{0}+\dfrac{\lambda_{1}g(P^{'}Z)}{(\lambda_{3}-\lambda_{1})(\lambda_{2}-\lambda_{1})}\\
\overset{.}{z}_{1}&=\lambda_{2}z_{1}+\dfrac{\lambda_{2}g(P^{'}Z)}{(\lambda_{3}-\lambda_{2})(\lambda_{1}-\lambda_{2})}\\
\overset{.}{z}_{2}&=\lambda_{3}z_{2}+\dfrac{\lambda_{3}^{2}g(P^{'}Z)}{(\lambda_{3}-\lambda_{1})(\lambda_{3}-\lambda_{2})}
\end{cases}
\end{equation}
Now, we construct the following Lyapunov function:
\begin{equation}
V(Z)=\dfrac{1}{2}(\lambda_{2}\lambda_{3}\|z_{0}\|^{2}+\lambda_{1}\lambda_{3}\|z_{1}\|^{2}+\lambda_{1}\lambda_{2}\|z_{2}\|^{2}).
\end{equation}
It follows that the derivative of V along the trajectories of $(9)$, denoted by $\overset{.}{V}(Z)$, is given by
\begin{align*}
\overset{.}{V}(Z)=({\dfrac{\partial V(Z)}{\partial z_0}})^{T}\overset{.}z_0+({\dfrac{\partial V(Z)}{\partial z_1}})^{T}\overset{.}z_1+({\dfrac{\partial V(Z)}{\partial z_2}})^{T}\overset{.}z_2 \\
=\lambda_1\lambda_2\lambda_3 \bigg(\|Z\|^{2}+ \dfrac{(g(P^{'}Z))^T}{(\lambda_{3}-\lambda_{1})(\lambda_{2}-\lambda_{1})} z_{0}+\\
\dfrac{(g(P^{'}Z))^T}{(\lambda_{3}-\lambda_{2})(\lambda_{1}-\lambda_{2})} z_{1}+\lambda_{3} \dfrac{(g(P^{'}Z))^T}{(\lambda_{3}-\lambda_{1})(\lambda_{3}-\lambda_{2})} z_{2}\bigg).\end{align*}

Hence by using Cauchy inequality and the Lipschitz property of $g$, we get
\begin{equation}
\begin{split}
\overset{.}{V}(Z)\le \lambda_1\lambda_2\lambda_3(1-L\|P^{'}\| \phi(\Lambda))\|Z\|^{2}
\end{split}
\end{equation}
where $\phi(\Lambda)$ is defined in Remark $1$.

 Next, we proceed to estimate the upper bound of $\|P^{'}\|$, where the matrix norm $\|.\|$ is the operator norm induced by the Euclidean norm, i.e. $\|P^{'}\|=\underset{\|w\|= 1}{\sup}{\|P^{'}w\|}$.

For any $w=(w_{1}^{T},w_{2}^{T},w_{3}^{T})^{T}\in \mathbb{R}^{3n}$  with $\|w\|=1$ where $w_{i}\in \mathbb{R}^{n}$, then by the
definition of $P^{'}$, we have $$P^{'}w=\begin{bmatrix}w_{1}+w_{2}+\frac{1}{\lambda_{3}}w_{3}\\\lambda_{1}w_{1}+\lambda_{2}w_{2}+w_{3}\end{bmatrix}.$$
By using Minkowski inequality and Cauchy inequality, we obtain
\begin{align*}
\|P^{'}w\|^{2}&= (\|w_{1}+w_{2}+\frac{1}{\lambda_{3}}w_{3}\|^{2}+\|\lambda_{1}w_{1}+\lambda_{2}w_{2}+w_{3}\|^{2})\\
&\le (3+\frac{1}{\lambda_{3}^{2}}+\lambda_{1}^{2}+\lambda_{2}^{2})(\|w_{1}\|^{2}+\|w_{2}\|^{2}+\|w_{3}\|^{2})\\
&=(3+\lambda_{1}^{2}+\lambda_{2}^{2}+\frac{1}{\lambda_{3}^{2}})
 .\end{align*}
Hence $\|P^{'}\|\le \sqrt{3+\lambda_{1}^{2}+\lambda_{2}^{2}+\frac{1}{\lambda_{3}^{2}}}$, which is $h(\Lambda)$ in Remark $1$ by definition. Consequently,
 from $(11)$ and the above fact, we know that
\begin{equation}
\overset{.}{V}(Z) \le \lambda_{1}\lambda_{2}\lambda_{3}(1-L\phi(\Lambda)h(\Lambda))\|Z\|^{2}
\end{equation}
for any $Z\in \mathbb{R}^{3n}$.

Notice that if the parameters $(k_{p},k_{i},k_{d})$ are chosen from $\Omega_{K}$, then the corresponding parameters $(\lambda_{1},\lambda_{2},\lambda_{3})$ should belong to $\Omega_{\Lambda}$, where $$\Omega_{\Lambda}=\{\Lambda \big | L\phi(\Lambda)h(\Lambda)<1; \lambda_{i}<0,i=1,2,3; \lambda_{i}\neq\lambda_{j},i\neq j\}.$$
Consequently, the right hand side of $(12)$ is a negative definite quadratic form of $Z$. Therefore,
for any $Z(0)\in \mathbb{R}^{3n}$, we have $\lim_{t\to\infty} Z(t)=0$ with an exponential rate of convergence from the Lyapunov stability theory, which in turn gives $$\lim_{t\to\infty} Y(t)=\lim_{t\to\infty} PZ(t)=0.$$
Recall that $Y(t)=(y_{0}(t),y_{1}(t),y_{2}(t))$ and $y_{1}(t)=x_{1}(t)-y^{*}, y_{2}(t)=x_{2}(t)$, thus we have $\lim_{t \to \infty} x_{1}(t)=y^{*}$ and $\lim_{t \to \infty} x_{2}(t)=0$ with an exponential rate of convergence, for any initial values.

Finally, to complete the proof of Theorem $1$, we show that $\Omega_{\Lambda}$ is an unbounded open set in $\mathbb{R}^{3}$. Let us choose two distinct negative numbers $\lambda_{1},\lambda_{2}$ arbitrarily, it is easy to see that $\phi(\Lambda)$ tends to $0$ and $h(\Lambda)$ keeps bounded as $\lambda_{3}$ tends to $-\infty$. Hence $\Omega_{\Lambda}$ is  nonempty and unbounded. Notice that $\Omega_{\Lambda}$ is the preimage of an open set of a continuous function, therefore it is open.

 Furthermore, from the relationship,
\begin{equation}\begin{cases}
k_{p}&=-(\lambda_{1}\lambda_{2}+\lambda_{1}\lambda_{3}+\lambda_{2}\lambda_{3})\\
k_{i}&=\lambda_{1}\lambda_{2}\lambda_{3}\\
k_{d}&=\lambda_{1}+\lambda_{2}+\lambda_{3}
\end{cases}.
\end{equation}

 We claim that $\Omega_{K}$ is also an open set in $\mathbb{R}^{3}$ since the Jacobian matrix of the mapping defined by $(13)$ is nonsingular at every point $\Lambda\in\Omega_{\Lambda}$, i.e.
\begin{align*}
\det\begin{bmatrix}
\frac{\partial k_{p}}{\partial{\lambda_{1}}}&\frac{\partial k_{p}}{\partial{\lambda_{2}}}&\frac{\partial k_{p}}{\partial{\lambda_{3}}}\\
\frac{\partial k_{i}}{\partial{\lambda_{1}}}&\frac{\partial k_{i}}{\partial{\lambda_{2}}}&\frac{\partial k_{i}}{\partial{\lambda_{3}}}\\
\frac{\partial k_{d}}{\partial{\lambda_{1}}}&\frac{\partial k_{d}}{\partial{\lambda_{2}}}&\frac{\partial k_{d}}{\partial{\lambda_{3}}}\\\end{bmatrix}
\end{align*}
\begin{align*}
&=\det \begin{bmatrix}
-(\lambda_{2}+\lambda_{3})&-(\lambda_{1}+\lambda_{3})&-(\lambda_{1}+\lambda_{2})\\
\lambda_{2}\lambda_{3}&\lambda_{1}\lambda_{3}&\lambda_{1}\lambda_{2}\\
1&1&1
\end{bmatrix}\\&=(\lambda_{1}-\lambda_{2})(\lambda_{1}-\lambda_{3})(\lambda_{3}-\lambda_{2})\neq0.
\end{align*}

Obviously, $\Omega_{K}$ is unbounded.

This completes the proof of the theorem.

\textbf{Proof of Corollary 1:}
From the previous analysis and the relationship $(13)$, it is sufficient to show that $(-\alpha,-(1+\alpha),-\beta)\in \Omega_{\Lambda} $ for  $0<\alpha<\frac{1}{4}$ and $\beta>\max\{5L,5\}$.

Notice that \begin{align*}\phi(-\alpha,-(1+\alpha),-\beta)&\le\phi(-\frac{1}{4},-\frac{5}{4},-\beta)\\ &< \sqrt {\frac{5}{\beta^{2}}}\end{align*} and that $$h(-\alpha,-(1+\alpha),-\beta)\le \phi(-\frac{1}{4},-\frac{5}{4},-5)< \sqrt5$$ whenever $0<\alpha<\frac{1}{4}$ and $\beta>5$.
Hence if $0<\alpha<\frac{1}{4}$ and $\beta>\max\{5L,5\}$, then we get $$\phi(-\alpha,-(1+\alpha),-\beta)h(-\alpha,-(1+\alpha),-\beta)<\frac{1}{L}.$$
This means that $\Omega_{K}^{'}\subset\Omega_{K}$.

\textbf{Proof of Proposition 2:}
Let $\epsilon >0$,  $y^{*}\in \mathbb{R}$ and $k_{i},k_{p},k_{d}$ be arbitrarily given.
Denote $y_{0}(t)=\int_{0}^{t} e(s)ds$, $y_{1}(t)=e(t)$, $y_{2}(t)=\overset{.}{e}(t)$, $y(t)=(y_{0}(t),y_{1}(t),y_{2}(t))^{T}$, then $(4)$ is equivalent to
\begin{equation}
\begin{cases}
\overset{.}{y}_{0}=y_{1}\\
\overset{.}{y}_{1}=y_{2}\\
\overset{.}{y}_{2}=((y_{1}+y^{*})^{2}+y_{2}^{2})^{\dfrac{1+\epsilon}{2}} +k_{i}y_{0}+k_{p}y_{1}+k_{d}y_{2}
\end{cases}
\end{equation}

For $L>1$, we define an unbounded closed set $C_{L}= \{y=(y_{0},y_{1},y_{2})\in \mathbb{R}^{3} \big| y_{2}-1 \ge y_{1} \ge y_{0}+L \ge L\}$. The shape of the set $C_{L}$ is a closed cone with the point $(0,L,L+1)$ being the vertex of this cone. Its boundary consists of three angular domain. They are

 $$\begin{cases}S_{1}=\{y \big| y_{2}-1 \ge y_{1} \ge y_{0}+L = L\} \\
S_{2}=\{y \big| y_{2}-1 \ge y_{1} = y_{0}+L \ge L\} \\
S_{3}=\{y \big| y_{2}-1 = y_{1} \ge y_{0}+L \ge L\}
 \end{cases}.$$

Define three vectors $v_{1}=(1,0,0)$,  $v_{2}=(-1,1,0)$ and  $v_{3}=(0,-1,1)$, which are the normal vectors of $S_{1}, S_{2}, S_{3}$, pointing to the inner side of the cone respectively.

We claim that, for large $L$, the following two inequalities
\begin{equation}
((y_{1}+y^{*})^{2}+y_{2}^{2})^{\dfrac{1+\epsilon}{2}} +k_{i}y_{0}+k_{p}y_{1}+k_{d}y_{2} \ge y_{2}+\dfrac{y_{2}^{1+\epsilon}}{2}
\end{equation}

and
\begin{equation}
((y_{1}+y^{*})^{2}+y_{2}^{2})^{\dfrac{1+\epsilon}{2}} +k_{i}y_{0}+k_{p}y_{1}+k_{d}y_{2} \le (3y_{2}^{2})^{\dfrac{1+\epsilon}{2}}
\end{equation}
holds for any $(y_{0},y_{1},y_{2}) \in C_{L}$.

This is easy to verify, because
$((y_{1}+y^{*})^{2}+y_{2}^{2})^{\dfrac{1+\epsilon}{2}} +k_{i}y_{0}+k_{p}y_{1}+(k_{d}-1)y_{2}-\dfrac{y_{2}^{1+\epsilon}}{2} \ge
\dfrac{y_{2}^{1+\epsilon}}{2}-(|k_{i}|+|k_{p}|+|k_{d}|+1)y_{2}$ for any $y \in C_{L}$. Hence we can take $L_{1}=\epsilon^{-1}(\log2(|k_{i}|+|k_{p}|+|k_{d}|+1))$ such that $(15)$ holds for $y \in C_{L_{1}}$. Similarly, we can choose $L_{2}$ to ensure $(16)$ holds for $y \in C_{L_{2}}$. Then  $L=\max \{L_1,L_2\}$ satisfies our claim.

Next, we consider the vector field defined by the equation $(14)$, we claim that: $C_{L}$ is an invariant set for $(14)$. This means that for any initial value $y(0)=(y_{0}(0),y_{1}(0),y_{2}(0))$ lies in $C_{L}$, the positive trajectory of  $y(0)$, which is defined as $\{y(t)\in \mathbb{R}^{3} \big| 0 \le t < a \}$, is contained in $C_{L}$ where $[0,a)$ is the maximal existence interval of the solution. It may equal to $\infty$ or perhaps a finite number. (We will prove shortly that, only  the latter case can happen.)

To prove the claim, it is sufficient to work with those initial values which lie in the boundary of $C_{L}$ because $(14)$ is an autonomous system.

Case $1$: For initial values lie in $S_{1}$, we have $v_{1}.\overset{.}{y}=y_{1}\ge L>0$.

Case $2$: For initial points lie in $S_{2}$, we have $v_{2}.\overset{.}{y}=y_{1}.-1+y_{2}.1=y_{2}-y_{1}=1>0$.

Case $3$: For initial points lie in $S_{3}$, we have
$$v_{3}.\overset{.}{y}=((y_{1}+y^{*})^{2}+y_{2}^{2})^{\dfrac{1+\epsilon}{2}} +k_{i}y_{0}+k_{p}y_{1}+(k_{d}-1)y_{2}>0$$ from $(15)$.

Hence for any initial value lies in $C_{L}$, the trajectory will not escape from the set $C_{L}$. On the other hand, as long as the trajectory lies in $C_{L}$, then we have $\overset{.}{e}(t)=y_{2}(t)\ge L+1$.

Now, take $x_{1}(0)=L+y^{*}$, $x_{2}(0)=L+1$, then from the relationship $y_{1}(t)=e(t)=x_{1}(t)-y^{*}, y_{2}(t)=\overset{.}{e}(t)=x_{2}(t)$ and $y_{0}(t)=\int_{0}^{t} e(s)ds$, we get $y(0)=(y_{0}(0), y_{1}(0),y_{2}(0))=(0,L,L+1)\in C_{L}$, thus from the above analysis, we have $y(t)\in C_{L}$ for any $t$ in the interval where  $(14)$ exists.

Note that $y_{2}(0)=L+1$, from $(15)$, we can get  $\overset {.}{y}_{2}\ge \dfrac{y_{2}^{1+\epsilon}}{2}$, then $y_{2}(t)\ge ((L+1)^{-\epsilon}-\dfrac{ t \epsilon}{2})^{- \dfrac{1}{\epsilon}}$ in the interval where $(14)$ exists from the comparison theorem. Denote $[0,a)$ be the maximal existence interval of the solution, we have $a\le \dfrac{2}{\epsilon(L+1)^{\epsilon}} < \infty$.
At the same time,  note that $\overset{.}{e}=y_{2}(t)\ge L+1$, we get $e(t)\ge  L+(L+1)t$ when $t\in [0,a)$.  This proves the first half of Proposition $2$.

 If $\epsilon\le 1$ and the initial value is $(0,L,L+1)$,  we claim that $\lim_{t \to a} y_{2}(t)=\infty$. Otherwise, we get $\lim_{t \to a} y_{1}(t)< \infty$ and hence $\lim_{t \to a}y_{0}(t)< \infty$ from the finiteness of $a$, which is a contradiction.

Now we are in a position to prove the next half of this proposition. From $(14)$ and $(16)$, we get
\begin{align*}
\dfrac{dy_{1}}{dy_{2}}&=y_{2}/(((y_{1}+y^{*})^{2}+y_{2}^{2})^{\dfrac{1+\epsilon}{2}} +k_{i}y_{0}+k_{p}y_{1}+k_{d}y_{2})\\
&\ge \dfrac{c_{\epsilon}y_{2}}{y_{2}^{1+\epsilon}}=\dfrac{c_{\epsilon}}{y_{2}^{\epsilon}}.
\end{align*}
Hence
\begin{align*}
\lim_{t \to a}e(t)=\lim_{t \to a}y_{1}(t)&=y_{1}(0)+\int_{L+1}^{\infty} \dfrac{dy_{1}}{dy_{2}} dy_{2}\\&\ge y_{1}(0)+\int_{L+1}^{\infty} \dfrac{c_{\epsilon}}{y_{2}^{\epsilon}} dy_{2}\\&= \infty
.\end{align*}
This completes our proof.

\textbf{Proof of Proposition 3:} Denote

 $y(t)=(y_{0}(t),\cdots,y_{3}(t))^{T}=(\int_{0}^{t} e(s)ds,e(t),\cdots,e^{(2)}(t))$. Let $f(x_{1},x_{2},x_3)=c  x_3$ where $0 < c \le L$ %is a number undetermined
 , then $f \in\mathscr{G}_L $, and $(5)$ turns to be
\begin{equation}
\begin{cases}
\overset{.}{y}_{0}&=y_1\\
\overset{.}{y}_1&=y_2\\
\overset{.}{y}_2&=y_3\\
\overset{.}{y}_3&=k_i y_{0} +k_p y_{1}+k_d y_{2}+cy_{3}.
\end{cases}
\end{equation}
Denote $A=\begin{bmatrix}
0&1&0&0&\\
0&0&1&0&\\
0&0&0&1&\\
k_{i}&k_{p}&k_{d}&c
\end{bmatrix}$,
then $(17)$ becomes
\begin{equation}
\overset{.}{y}=Ay.
\end{equation}
It is easy to see that the characteristic polynomial of $A$ equals to $\lambda^{4}-c\lambda^{3}-k_{d}\lambda^{2}-k_{p}\lambda-k_{i}$. Denote $\lambda_{0},\cdots,\lambda_{3}$ be the $4$  complex eigenvalues of $A$. Then
\begin{equation}
\lambda_{0}+\cdots+\lambda_{3}=c
\end{equation}
for any choice of $k_{p},k_{i},{k_{d}}$.
Take the real parts of $(19)$, we see that there exists at least one eigenvalue whose real part is strictly positive.
Hence $A$ is not a Hurwitz matrix for any $k_{p},k_{i},k_{d}$.

Now, we prove Proposition $3$ in the case $k_{i}\neq 0$ first. By Lemma A in the Appendix, we can take $c$ to make the matrix $A$ be similar to
 $$J=\begin{bmatrix}
 \lambda_{0}&0&0&0\\
0&\lambda_{1}&0&0\\
0&0&\lambda_{2}&0\\
0&0&0&\lambda_{3}
\end{bmatrix}.$$

Let us denote $$P=\begin{bmatrix}
 1&1&1&1\\
\lambda_{0}&\lambda_{1}&\lambda_{2}&\lambda_{3}\\
\lambda_{0}^{2}&\lambda_{1}^{2}&\lambda_{2}^{2}&\lambda_{3}^{2}\\
 \lambda_{0}^{3}&\lambda_{1}^{3}&\lambda_{2}^{3}&\lambda_{3}^{3}\\
\end{bmatrix},$$  then $A=PJP^{-1}$.
Without loss of generality, we suppose that $\Re(\lambda_{0}) \le \Re(\lambda_{1}) \le \Re(\lambda_{2}) \le \Re(\lambda_{3})$, where $\Re(.)$ denote the real part of a complex number.

Obviously, the solution of $(17)$ can be expressed by the formula
\begin{equation}
 y(t)=e^{At}y(0).
\end{equation}

 From this, we get $y(t)=e^{At}y(0)=Pe^{Jt}P^{-1}y(0)$.

Denote $z(0)=P^{-1}y(0)$, and choose $$y(0)=P\begin{bmatrix}
0\\0\\-1\\1
\end{bmatrix}$$ such that $z(0)=(0,0,-1,1)^{T}$.  By simple calculations, we know that $$y(t)=Pe^{Jt}\begin{bmatrix}
0\\0\\-1\\1
\end{bmatrix}=\begin{bmatrix}
e^{\lambda_{0}t}&e^{\lambda_{1}t}&e^{\lambda_{2}t}&e^{\lambda_{3}t}\\
\lambda_{0}e^{\lambda_{0}t}&\lambda_{1}e^{\lambda_{1}t}&\lambda_{2}e^{\lambda_{2}t}&\lambda_{3}e^{\lambda_{3}t}\\
\lambda_{0}^{2}e^{\lambda_{0}t}&\lambda_{1}^{2}e^{\lambda_{1}t}&\lambda_{2}^{2}e^{\lambda_{2}t}&\lambda_{3}^{2}e^{\lambda_{3}t}\\
\lambda_{0}^{3}e^{\lambda_{0}t}&\lambda_{1}^{3}e^{\lambda_{1}t}&\lambda_{2}^{3}e^{\lambda_{2}t}&\lambda_{3}^{3}e^{\lambda_{3}t}\\
\end{bmatrix} \begin{bmatrix}
0\\0\\-1\\1
\end{bmatrix}.$$
Hence $y_{1}(t)=\lambda_{3}e^{\lambda_{3}t}-\lambda_{2}e^{\lambda_{2}t}$.

We now show that $e(t)$ is unbounded by considering three cases separately.

Case 1: If $\Re(\lambda_{3})> \Re(\lambda_{2})$, then it is obviously that $$\lim_{t \to \infty} |y_{1}(t)|=\infty.$$

Case 2: If $\Re(\lambda_{3})= \Re(\lambda_{2})=a$ and $|\lambda_{2}| \neq |\lambda_{3}| $ , then  $$|y_{1}(t)|=|\lambda_{3}e^{\lambda_{3}t}-\lambda_{2}e^{\lambda_{2}t}| \ge e^{at}| |\lambda_{2}|-|\lambda_{3}| |,$$
we also get $\lim_{t \to \infty} |y_{1}(t)|=\infty$ in this case.

Case 3: If $\Re(\lambda_{3})= \Re(\lambda_{2})=a$ and $|\lambda_{2}| = |\lambda_{3}|$,  then $\lambda_{2} ,\lambda_{3}$ are conjugate complex numbers.
Let $\lambda_{3}=a+bi,\lambda_{2}=a-bi$, then $\lambda_{3}e^{\lambda_{3}t}-\lambda_{2}e^{\lambda_{2}t}=e^{(a-bi)t}(\lambda_{3}e^{2bit}-\lambda_{2}).$ It is unbounded, too.

Thus, in any case, we get the unboundedness of $e(t)$. Hence we have proved our proposition in the case $k_{i}\neq 0$.

Next, if $k_{i}=0$, then $(17)$ reduces to
\begin{equation}
\begin{cases}
\overset{.}{y}_1&=y_2\\
\overset{.}{y}_2&=y_3\\
\overset{.}{y}_3&=k_p y_{1}+k_d y_{2}+cy_{3}.
\end{cases}
\end{equation}
The corresponding conclusion can be proven analogously.
Details are omitted.

\section{Conclusions}
  In this paper, we have presented a theoretical investigation on the   capability and limitations of the
 widely used classical PID controller for a basic class of nonlinear  uncertain dynamical systems,  found a three dimensional
manifold (in Theorem $1$) within which the three controller parameters can be taken arbitrarily to design a globally stabilizing PID controller.  To the best of
our knowledge,  this appears to be the first of such kind of results in the literature on  PID  control. Of course,  many interesting problems
still remain open. It would be interesting to give a comparison with the active disturbance rejection control (ADRC) in \cite{Han2009}.  It would also be interesting to know whether or not the  PID controller (and its extensions) can be applied to  nonlinear system
structures more complicated than, for example, those  treated by the back-stepping design (see \cite{Krstic1995}). Furthermore, it is desirable to consider time-delay systems  and
sampled-data PID controllers under a prescribed sampling rate,  within the related boundaries established for the maximum capability of the general feedback mechanism
( cf. e.g. \cite{Xie2000}, \cite{Ren2014}). These belongs to further investigation.

\section{appendix}

 Lemma A: For any $L >0$, and any $k_{i}\neq0$ and any $k_{p}, k_{d} \in \mathbb{R}$, there exist $c\in (0,L]$, such that the equation $\lambda^{4}-c\lambda^{3}-k_{d}\lambda^{2}-k_{p}\lambda-k_{i}=0$ has four distinct roots in the complex plane.

\textbf{Proof of lemma A:} Denote $g(\lambda)=\lambda^{4}-c\lambda^{3}-k_{d}\lambda^{2}-k_{p}\lambda-k_{i}$,
$A_{c}=\{\lambda:\lambda^{4}-c\lambda^{3}-k_{d}\lambda^{2}-k_{p}\lambda-k_{i}=0\}$. Then it is easy to verify that $A_{c_{1}} $ and $A_{c_{2}} $ are disjoint when $c_{1}\neq c_{2}$.

It is well known that, the sufficient and necessary condition of a polynomial $g(\lambda)$ has multiple roots is that $g(\lambda)$ has at least one common zeros with its derivative $g^{'}(\lambda)$,
$$g^{'}(\lambda)=4\lambda^{3}-3c\lambda^{2}-2k_{d}\lambda-k_{p}.$$
If $w$ is a multiple root of $g$, then
\begin{equation}
\begin{cases}
w^{4}-cw^{3}-k_{d}w^{2}-k_{p}w-k_{i}=0\\
4w^{3}-3cw^{2}-2k_{d}w-k_{p}=0
\end{cases}
.\end{equation}
From $(22)$, we get
\begin{equation}
w^{4}+k_{d}w^{2}+2k_{p}w+3k_{i}=0.
\end{equation}Denote $R=\{w:w^{4}+k_{d}w^{2}+2k_{p}w+3k_{i}=0\}$.
From the above analysis, we get $w\in A_{c}\cap R$ iff $w$ is a multiple roots of $g(\lambda)$.
There are at most 4 elements in the set $R$ and $R$ is independent of $c$.
Hence we can choose $c\in (0,L]$ such that $A_{c}\cap R=\varnothing$. Such $c$ will make $g(\lambda)$ to have no multiple roots. The proof of the lemma is complete.

\bibliography{ifacconf}

\begin{thebibliography}{17}
\providecommand{\natexlab}[1]{#1}
\providecommand{\url}[1]{\texttt{#1}}
\providecommand{\urlprefix}{URL }
\expandafter\ifx\csname urlstyle\endcsname\relax
  \providecommand{\doi}[1]{doi:\discretionary{}{}{}#1}\else
  \providecommand{\doi}{doi:\discretionary{}{}{}\begingroup
  \urlstyle{rm}\Url}\fi

\bibitem[{{\AA}str{\"o}m and H{\"a}gglund(1995)}]{Astrom1995}
{\AA}str{\"o}m, K.J. and H{\"a}gglund, T. (1995).
\newblock \emph{\protect{PID Controllers: Theory, Design and Tuning}}.
\newblock Instrument society of America.

\bibitem[{{\AA}str{\"o}m and H{\"a}gglund(2006)}]{Astrom2006}
{\AA}str{\"o}m, K.J. and H{\"a}gglund, T. (2006).
\newblock \emph{\protect{Advanced PID Control}}.
\newblock ISA-The Instrumentation, Systems and Automation Society.

\bibitem[{Blanchini et~al.(2004)Blanchini, Lepschy, Miani, and
  Viaro}]{Blanchini2004}
Blanchini, F., Lepschy, A., Miani, S., and Viaro, U. (2004).
\newblock \protect{Characterization of PID and lead/lag compensators satisfying
  given $H_{\infty}$ specifications}.
\newblock \emph{IEEE Transactions on Automatic Control}, 49(5), 736--740.

\bibitem[{Guo(2011)}]{Guo2011}
Guo, L. (2011).
\newblock \protect{Some perspectives on the development of control theory}.
\newblock \emph{J.Sys. Sci and Math. Scis.}, 31(9), 1014--1018.

\bibitem[{Guo(2014)}]{Guo2014}
Guo, L. (2014).
\newblock \protect{How much uncertainty can feedback mechanism deal with?}
\newblock \emph{Plenary Lecture at the 19th IFAC World Conngress, Cape Town.
  July, 2014 www.ifac2014.org/assets/pdf/plenary/Guo.pdf}.

\bibitem[{Han(2009)}]{Han2009}
Han, J. (2009).
\newblock \protect{From PID to active disturbance rejection control}.
\newblock \emph{IEEE transactions on Industrial Electronics}, 56(3), 900--906.

\bibitem[{Hara et~al.(2006)Hara, Iwasaki, and Shiokata}]{Hara2006}
Hara, S., Iwasaki, T., and Shiokata, D. (2006).
\newblock \protect{Robust PID control using generalized KYP synthesis: Direct
  open-loop shaping in multiple frequency ranges}.
\newblock \emph{IEEE control systems}, 26(1), 80--91.

\bibitem[{Ho et~al.(1997)Ho, Datta, and Bhattacharyya}]{Ho1997}
Ho, M.T., Datta, A., and Bhattacharyya, S. (1997).
\newblock \protect{A linear programming characterization of all stabilizing PID
  controllers}.
\newblock In \emph{American Control Conference, 1997}, volume~6, 3922--3928.

\bibitem[{Ho and Lin(2003)}]{Ho2003}
Ho, M.T. and Lin, C.Y. (2003).
\newblock \protect{PID controller design for robust performance}.
\newblock \emph{IEEE Transactions on Automatic Control}, 48(8), 1404--1409.

\bibitem[{Keel and Bhattacharyya(2008)}]{Keel2008}
Keel, L.H. and Bhattacharyya, S.P. (2008).
\newblock Controller synthesis free of analytical models: Three term
  controllers.
\newblock \emph{IEEE Transactions on Automatic Control}, 53(6), 1353--1369.

\bibitem[{Killingsworth and Krstic(2006)}]{Killingsworth2006}
Killingsworth, N.J. and Krstic, M. (2006).
\newblock \protect{PID tuning using extremum seeking: online, model-free
  performance optimization}.
\newblock \emph{IEEE control systems}, 26(1), 70--79.

\bibitem[{Krstic et~al.(1995)Krstic, Kanellakopoulos, and
  Kokotovic}]{Krstic1995}
Krstic, M., Kanellakopoulos, I., and Kokotovic, P.V. (1995).
\newblock \emph{Nonlinear and Adaptive Control Design}.
\newblock Wiley.

\bibitem[{Ren et~al.(2014)Ren, Cheng, and Guo}]{Ren2014}
Ren, J., Cheng, Z., and Guo, L. (2014).
\newblock \protect{Further results on limitations of sampled-data feedback}.
\newblock \emph{Journal of Systems Science and Complexity}, 27(5), 817--835.

\bibitem[{Silva et~al.(2007)Silva, Datta, and Bhattacharyya}]{silva2007pid}
Silva, G.J., Datta, A., and Bhattacharyya, S.P. (2007).
\newblock \emph{\protect{PID Controllers for Time-Delay Systems}}.
\newblock Springer Science \& Business Media.

\bibitem[{S{\"o}ylemez et~al.(2003)S{\"o}ylemez, Munro, and
  Baki}]{Soylemez2003}
S{\"o}ylemez, M.T., Munro, N., and Baki, H. (2003).
\newblock \protect{Fast calculation of stabilizing PID controllers}.
\newblock \emph{Automatica}, 39(1), 121--126.

\bibitem[{Xie and Guo(2000)}]{Xie2000}
Xie, L.L. and Guo, L. (2000).
\newblock \protect{How much uncertainty can be dealt with by feedback?}
\newblock \emph{IEEE Transactions on Automatic Control}, 45(12), 2203--2217.

\bibitem[{Ziegler and Nichols(1942)}]{ziegler1942}
Ziegler, J.G. and Nichols, N.B. (1942).
\newblock Optimum settings for automatic controllers.
\newblock \emph{Trans. ASME}, 64(11).

\end{thebibliography}

\end{document}